\documentclass[letterpaper, 10pt, conference]{ieeeconf}
\IEEEoverridecommandlockouts                              
\usepackage[utf8]{inputenc}
\usepackage{epsfig}
\usepackage{amsmath}

\usepackage{amsthm}
\usepackage{amsfonts}
\usepackage{eso-pic}
\usepackage{textcomp}
\usepackage{comment}
\usepackage{graphicx}
\usepackage{amssymb}
\usepackage{gensymb}
\usepackage[font=small,skip=0pt]{caption}
\usepackage[acronym]{glossaries}
\usepackage{algorithm}
\usepackage{algpseudocode}
\usepackage{xcolor}
\graphicspath{{Figures/}}
\usepackage{draftwatermark}

\makeatletter
\let\NAT@parse\undefined
\makeatother
\usepackage{hyperref}
\hypersetup{
colorlinks=true, 
breaklinks=true, 
urlcolor= black, 
citecolor=black,
linkcolor= blue, 
bookmarksopen=true,
pdftitle={Autonomous driving using an optimized neural network based adaptive LPV-MPC controller}, 
pdfauthor={Y. Kebbati},
}

\title{\LARGE \bf
Autonomous driving using an optimized neural network based adaptive LPV-MPC controller
}

\author{Yassine Kebbati$^{1*}$, Naima Ait-Oufroukh$^{1}$, Vicenç Puig$^{2}$, Vincent Vigneron$^{1}$ and Dalil Ichalal$^{1}$
\thanks{$^*$Corresponding author: \href{mailto:yassine.kebbati@univ-evry.fr}{yassine.kebbati@univ-evry.fr}}
\thanks{$^1$IBISC-EA4526, univ Evry, université Paris-Saclay, France}
\thanks{$^2$CS2AC, Polytechnic University of Catalonia, Barcelona, Spain}
}

\makeatletter
\newcommand*{\rom}[1]{\expandafter\@slowromancap\romannumeral #1@}
\makeatother

\begin{document}

\maketitle

\thispagestyle{empty}
\pagestyle{empty}

\SetWatermarkText{Preprint}
\begin{abstract}
Driverless vehicles are complex systems operating in constantly changing environments. Automated driving is achieved by controlling the coupled longitudinal and lateral vehicle dynamics. Model predictive control is one of the most promising tools for this type of application due to its optimal performance and ability to handle constraints. This paper addresses autonomous driving with an adaptive linear parameter varying model predictive controller (LPV-MPC), which is adapted by a neural network and optimized by an improved Genetic Algorithm. The proposed controller is evaluated on a challenging track under variable wind disturbance. 

\end{abstract}

\section{Introduction}

Scientists and engineers in top automobile companies have been striving to achieve full driving autonomy by replacing human drivers with automatic control systems. Consequently, autonomous vehicles are emerging as a top technology aiming at improving traffic safety and enhancing mobility. The transition from manual to automatic driving is further accelerated by the huge leap in artificial intelligence and information processing technologies. Such a technological shift can also boost human productivity, in the sense that time spent on driving can be used for other productive tasks instead. Autonomous driving is a multidisciplinary field that involves sensing, perception, decision making and planning. In addition, accurate motion control is among the most important tasks for achieving full autonomy. The latter can be divided into longitudinal control in charge of speed tracking and lateral control which handles the steering task.

Literature research is divided into two categories; those who address the longitudinal and lateral control separately, and those who couple both tasks together. For instance, Xu \textit{et al.} developed an optimal preview controller for speed regulation that integrates road slope, speed profile and vehicle dynamics \cite{1}. Kebbati \textit{et al.} \cite{2} addressed speed control by designing a self-adaptive PID controller based on neural networks and genetic algorithms. On the other hand, Han \textit{et al.} \cite{3} designed an adaptive neural network PID controller for path tracking. They applied it to a second order vehicle model whose parameters were estimated by a forgetting factor least square algorithm. LPV $H_{\infty}$ was developed for high-speed driving and evasive maneuvers by Corno \textit{et al.} \cite{4}.  The design was based on the lateral error and look-ahead distance of the vehicle to ensure better robustness and account for actuator nonlinearities under low speeds. In \cite{5}, authors developed a model predictive controller (MPC) for lateral control and used fuzzy inference systems (FIS) to tune its weighting matrices. Guo \textit{et al}. \cite{6} developed a path-tracking MPC controller that considers the varying road conditions and small-angle assumptions as a form of measurable disturbance and solved the control problem using the differential evolution algorithm. Authors of \cite{7} worked on adaptive MPC for path tracking that optimizes the controller tuning with an improved PSO algorithm. Online controller adaptation was achieved by a lookup table approach, but this technique cannot account for all possible cases despite the good results that were obtained. However, the same authors improved their approach in \cite{8} by replacing the lookup table method with neural networks an adaptive neuro-fuzzy inference systems to generalize the adaptation beyond the data in the lookup table. Although significant tracking improvements were achieved, this approach still requires long offline optimizations. 

To address the mixed lateral and longitudinal control, Alcala \textit{et al.} developed in \cite{9} a solution for trajectory tracking. This control strategy was divided into a cascade control scheme, with an internal layer for controlling vehicle dynamics and an external one for vehicle kinematics. The same authors developed in \cite{10} a Takagi-Sugeno based MPC (TS-MPC) for autonomous driving. They used a data-driven approach to learn a Takagi-Sugeno representation of the vehicle dynamics, which was used in MPC with Moving Horizon Estimator (MHE) for achieving coupled longitudinal and lateral control. In \cite{11}, a coordinated lateral and longitudinal control using LPV-MPC for lateral control with PSO-PID for speed regulation is proposed that includes exponential weight to the MPC cost function to improve tracking performance. The use of nonlinear model predictive control (NMPC) was proposed in \cite{12} where the authors aimed at exploring a parking lot autonomously and performing the parking maneuver. Finally, an online learning MPC controller for autonomous racing by learning the model errors online using Gaussian process regression is proposed in \cite{13}.

This paper contributes to the above mentioned literature by developing an enhanced controller for the coupled longitudinal and lateral control. The contributions of this work are threefold: First, an adaptive LPV-MPC is developed for realizing autonomous driving. Second, an improved genetic algorithm is proposed for tuning the weighting matrices of the cost function to achieve optimal control actions. Third, a deep neural network is developed and trained to learn the tire lateral dynamics by predicting the cornering stiffness coefficients from measurable parameters such as velocities and accelerations. 

The article is organized as follows: Section \rom{2} presents the modeling of the vehicle dynamics. Section \rom{3} details the design of the LPV-MPC controller, the controller adaptation approach using machine learning and the controller optimization via the improved GA version. Evaluation results of the learning approach, optimization and control are presented and analyzed in section \rom{4}. Finally, section \rom{5} provides summarized conclusions and gives directions for future work. 

\section{Vehicle Modeling}
\begin{figure}[htb]
\centering
\includegraphics[width=7.75cm,height=5cm]{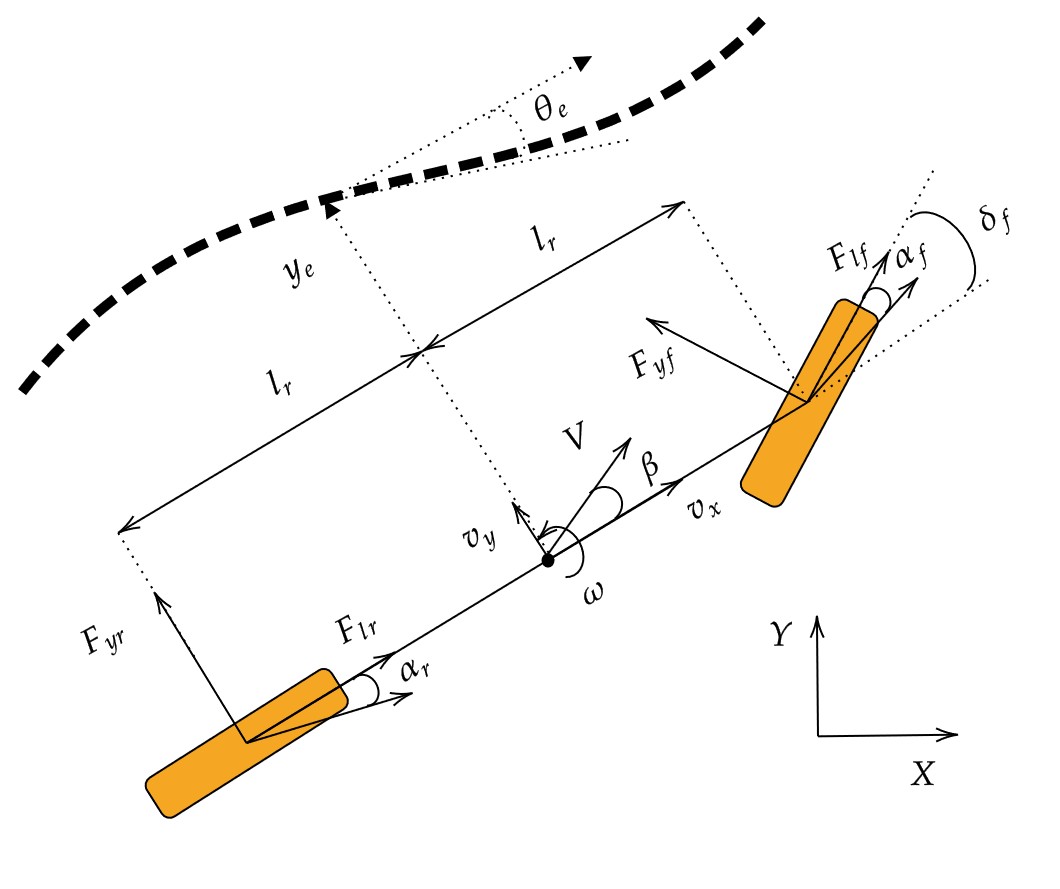}
\caption{Bicycle dynamic model with tracking error.}
\label{fig:1}
\end{figure}

The commonly used bicycle dynamics model \cite{14,15} is adopted in this paper as it is accurate for control design whilst being simple for real-time implementation. In this model, the two front wheels similar to the rear ones are lumped to form a single track representation (see Fig. \ref{fig:1}). The lateral dynamics are governed by the tire lateral forces which are a function of the slip angles. 
The full model accounts for longitudinal, lateral and yaw dynamics (\ref{eq1}) as well as tire forces (\ref{eq2}) and heading and lateral position errors (\ref{eq3}):
\begin{equation}
\label{eq1}
\left\{
    \begin{array}{ll}
        \dot{v}_x&= \alpha_x +\omega v_y - \frac{1}{m}(F_{yf}\sin{\delta}+F_d)\\
        \dot{v}_y&= \frac{1}{m}(F_{yf}\cos{\delta}+F_{yr}) -\omega v_x\\
        \dot{\omega}&= \frac{1}{I}(F_{yf}l_f\cos{\delta}-F_{yr}l_r)
    \end{array}
\right.
\end{equation}
\begin{equation}
\label{eq2}
\left\{
    \begin{array}{ll}
        F{yf}&= C_f \alpha_f\\
        F{yr}&= C_r \alpha_r\\
        F_d&= \mu mg + \frac{1}{2}\rho C_dAv_x^2
    \end{array}
\right.
\end{equation}
\begin{equation}
\label{eq3}
\left\{
    \begin{array}{ll}
        \dot{y}_e&= v_x\sin\theta_e + v_y \cos\theta_e\\
        \dot{\theta}_e&= \omega - \frac{v_x\cos\theta_e-v_y\sin\theta_e}{1-y_e k}
    \end{array}
\right.
\end{equation}
The parameters $v_x$, $v_y$ and $\omega$ represent the linear longitudinal and lateral velocities and yaw rate in the body frame. $F_{y(f,r)}$ are the tire lateral forces of the front and rear wheels respectively. $F_d$ is the total drag force where $C_d$ and $A$ represent the drag coefficient and the vehicle cross sectional area. $y_e$ and $\theta_e$ represent the lateral position error and the heading error where $k$ is the road curvature and $\epsilon$ is added to avoid singularities in the model. $I$ and $m$ are the inertia and the mass of the vehicle and $l_{(f,r)}$ are the distances between the vehicle centre of gravity and the front and rear wheel axles, respectively. $\alpha_x$ and $\delta$ are the acceleration and steering controls and $\mu$ and $g$ represent the friction coefficient and the gravity pull. Finally, $C_{(f,r)}$ represent the respective front and rear tire cornering stiffness coefficients, and $\alpha_{(f,r)}$ are the slip angles of the front and rear wheels which are respectively given by
\begin{equation}
\label{eq4}
\left\{
    \begin{array}{ll}
        \alpha_f&= \delta - \tan^{-1}{( \frac{v_y}{v_x+\epsilon}-\frac{l_f\omega}{v_x+\epsilon})}\\
        \alpha_r&= -\tan^{-1}{ (\frac{v_y}{v_x+\epsilon}+\frac{l_r\omega}{v_x+\epsilon})}
    \end{array}
\right.
\end{equation}
The full model can be summarized as a non-linear function of the state vector $(x)$, input vector $(u)$ and the road curvature $(k)$ as follows:
\begin{equation}
\label{eq5}
\dot{x} = f(x, u, k) 
\end{equation}
where $x = [v_x\ v_y\ \omega\ y_e\ \theta_e]^T$ and $u = [\delta\ \alpha_x]^T$.

\section{Controller Design}
The control task is achieved based on LPV-MPC approach \cite{14}, where the model presented in Section \rom{2} is reformulated in an LPV form. The model is then transformed into a state space representation where the state and control matrices are functions that depend on a scheduling vector of varying parameters. Embedding linear varying parameters into these matrices allows to capture the nonlinearities of the full model that provides a simple but accurate model for control design. Generally speaking, LPV systems are a class of linear systems whose parameters are functions of external or internal scheduling signals. The LPV state space formulation of the system is given by (\ref{eq6}) based on the scheduling vector $\psi=[\delta\ v_x\ v_y\ \theta_e\ y_e\ k ]^T$.
\begin{equation}
\label{eq6}
\dot{x} = A(\psi)x + B(\psi)u 
\end{equation}
According to \cite{14,16}, the state matrix $A(\psi)$ and control matrix $B(\psi)$ can be derived as the following:
\begin{equation}
\label{eq7}
A(\psi) = \left[{\begin{array}{ccccc}A_{11} & A_{12} & A_{13} & 0 & 0\\ 0 & A_{22} & A_{23} & 0 & 0\\0 & A_{32} & A_{33} & 0 & 0\\A_{41} & A_{42} & 0 & 0 & 0 \\ A_{51} & A_{52} & 1 & 0 & 0\end{array}}\right],
\end{equation}
\begin{equation}
\label{eq8}
B(\psi) = \left[\begin{array}{cc} B_{11} & 1 \\ B_{21} & 0 \\ B_{31} & 0\\ 0 & 0 \\ 0 & 0 \end{array}\right],
\end{equation}
where the terms are given by:\\
\noindent
$A_{11} = \frac{-\mu g}{v_x}- \frac{\rho C_d A v_x}{2m},\
A_{12} = \frac{C_f \sin{\delta}}{m v_x},\\
A_{13} = \frac{C_f l_f \sin{\delta}}{m v_x}+ v_y,\
A_{22} = -\frac{C_r + C_f \cos{\delta}}{m v_x},\\
A_{23} = -\frac{C_f l_f \cos{\delta} - C_r l_r}{m v_x}- v_x,\
A_{32} = -\frac{C_f l_f \cos{\delta + C_r l_r}}{I v_x},\\
A_{33} = -\frac{C_f l_f^2 \cos{\delta + C_r l_r^2}}{I v_x},\
A_{41} = \sin{\theta_e},\
A_{42} = \cos{\theta_e},\\
A_{45} = v_x,\
A_{51} = -\frac{k \cos{\theta_e}}{1-y_e k},\
A_{52} = \frac{k \cos{\theta_e}}{1 - y_e k},\\
B_{11} = - \frac{C_f \sin{\delta}}{m},\
B_{21} = -\frac{C_f \cos{\delta}}{m},\
B_{31} = -\frac{C_f l_f \cos{\delta}}{I}.$\\

Model predictive control uses the plant model to predict its behaviour over a prediction horizon $N_p$. Then, it generates an optimal control sequence by solving a constrained convex optimization problem. The receding horizon principle is then applied and only the first term of the optimal control sequence is used. The LPV model, discretized with $T_s$ sampling time, is used as the MPC prediction model. Hence, at each iteration the scheduling vector is used to instantiate the LPV model. The parameters of the scheduling vector can be obtained from sensors, planners or previous MPC predictions. In this regard, the MPC problem can be formulated as the following constrained quadratic cost function:
\begin{equation}
\label{eq9}
\begin{split}
\min_{\Delta U_k} \ J_k & = \sum_{i=0}^{N_p-1} \Big( (r_{k+i} - x_{k+i})^T Q (r_{k+i} - x_{k+i}) + \\
         & \Delta u_{k+i} R \Delta u_{k+i} \Big) + x_{k+N_p}^T Q x_{k+N_p}\\
s.t:\\ 
    & x_{k+i+1} = x_{k+i} + A(\psi_{k+i}) x_{k+i} + B(\psi_{k+i}) u_{k+i}dt \\
    & u_{k+i} = u_{k+i-1} + \Delta u_{k+i}\\
    & \Delta u_{min}  \leq \Delta u_k \leq \Delta u_{max}\\
    &  u_{min}  \leq u_k \leq \ u_{max}\\
    & x_{min}  \leq x_k \leq x_{max}
\end{split}
\end{equation}
The terms $x$, $u$ and $N_p$ are the previously defined state vector, control vector and prediction horizon. $Q \in \mathbb{R}^{5 \times 5}$ and $R \in \mathbb{R}^{2 \times 2}$ are semi-positive definite weighting matrices that penalise the states and the control effort. $r_{k+i}$ is the reference vector and the terms $[u_{min}, u_{max}]$, $[\Delta u_{min}, \Delta u_{max}]$ and $[x_{min}, x_{max}]$ are the upper and lower bounds on the control actions, control increments and states, respectively. 
\subsection{Controller adaptation with neural network} \label{sec:adaptation}
\begin{figure}[t]
\centering
\includegraphics[width=8cm,height=4.8cm]{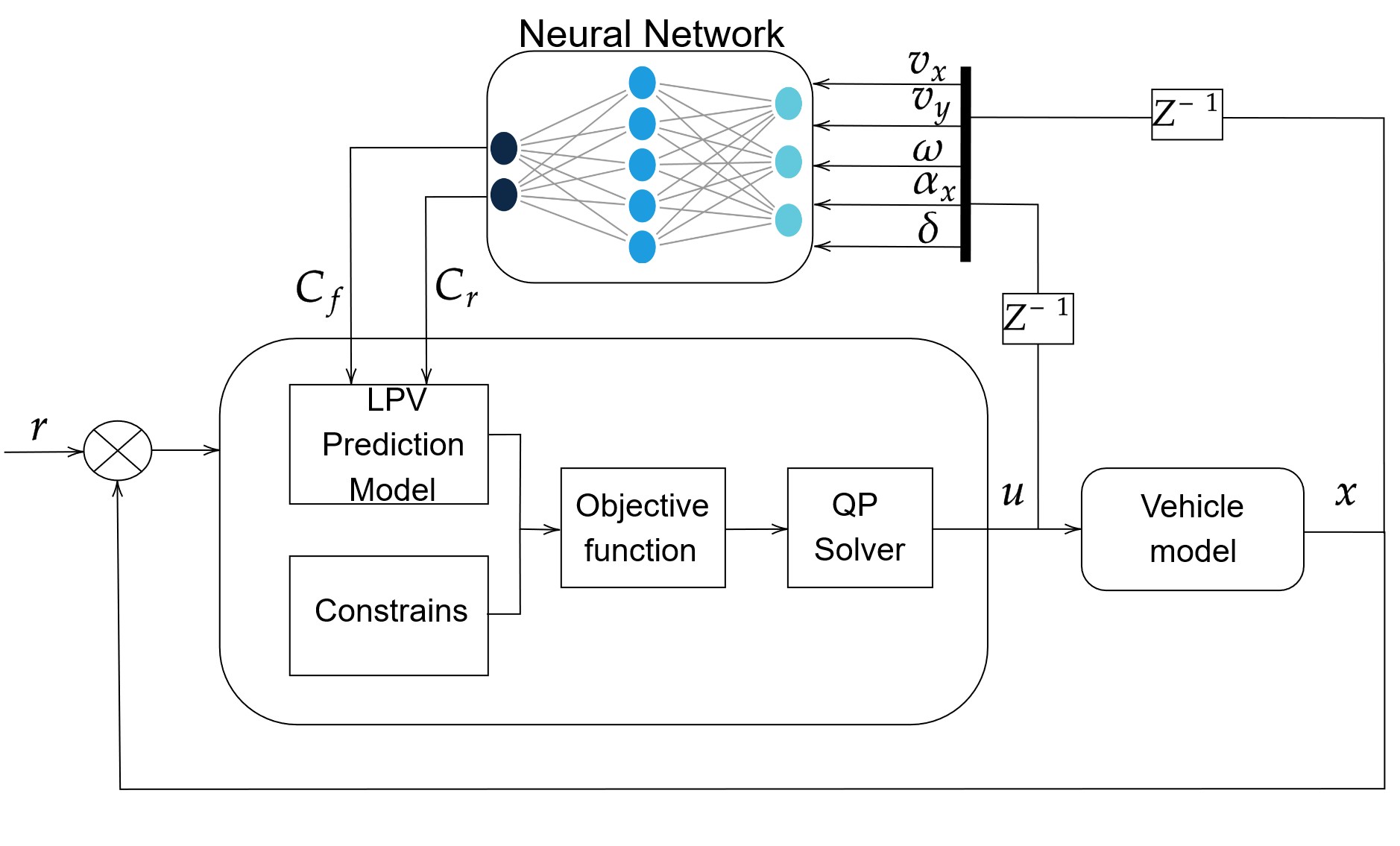}
\caption{Adaptive LPV-MPC approach.}
\label{fig:2}
\end{figure}

The majority of research in the literature consider a linearized tire model, where the relation between tire force and slip angle is governed by a constant called cornering stiffness coefficient. However, this is only valid for relatively small slip angles and not during fast and challenging maneuvers. To remedy this, we use machine learning to predict the cornering stiffness coefficient online from measurable parameters. We assume these parameters (namely, the longitudinal ($v_x$) and lateral ($v_y$) velocities, the steering angle ($\delta$), the acceleration ($\alpha_x$) and the yaw rate ($\omega$)) are enough to capture the tire dynamics. Using this approach (see Fig. \ref{fig:2}), the prediction model of the LPV-MPC is adapted online to improve the prediction capability and precision.        
\subsection{Controller tuning with improved GA Algorithm}
\begin{figure}[b]
\centering
\includegraphics[width=0.48\textwidth]{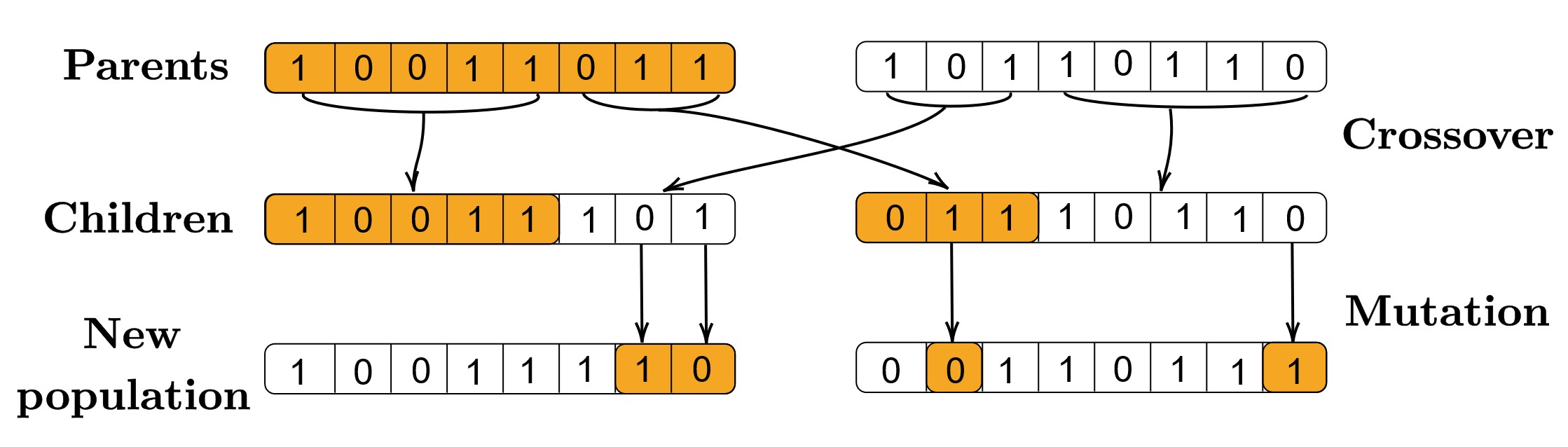}
\caption{Genetic operations}
\label{fig:3}
\end{figure}

Manually tuning the MPC controller can be time consuming and may not result in optimal performance. Thus, to optimize the designed LPV-MPC, we propose an improved Genetic Algorithm (GA) to tune the weighting matrices of the quadratic cost function. Genetic algorithms are a global search optimization technique based on biological evolution theory. They can easily find the optimum of an objective function even if is not continuous or differentiable \cite{17}. The algorithm initializes a population set of possible encoded solutions named chromosomes which are genetically enhanced over iterations and evaluated by a fitness function to determine their optimality.
The main operations in a GA are the selection, the crossover and the mutation processes, which control the search capability and the quality of the solutions. These operations consist of different functions which influence the algorithm at different degrees \cite{17}. The selection process chooses the best chromosomes to be enhanced by the crossover and mutation operations. The crossover seeks to produce high-quality solutions by mixing genetic data, while mutation introduces new genes to complement the crossover as illustrated in Fig. \ref{fig:3}. 

To improve the GA performance, researchers have essentially tried to improve the genetic operations. For instance, the most common selection operations are the roulette wheel (RWS) and tournament (TS) \cite{17}. Similarly, crossover versions include single/multi-point, uniform and shuffle crossover, while in mutation, one finds swap, inversion and random resetting. Details about these methods and GA can be found in \cite{19}. Rather than improving the operators themselves, we use a combination of RWS and TS selection operations to improve the selections of potential genes which is a critical phase of the GA. In particular, RWS strategy provides a higher chance for the good genes to be selected, which enhances exploitation and accelerates the convergence of the algorithm. However, RWS is mainly based on the fitness value, which makes it prone to premature convergence by possibly selecting the same dominant genes every time. On the other hand, the TS method allows to control the selection pressure. A smaller tournament size ensures more chances for weak genes to be selected unlike RWS. This feature retains diversity of the search pace and increases the possibility of converging to a global optimum at the expense of slower convergence. In our improved GA, both methods are used with random percentages at each iteration of the algorithm. The goal is to increase both convergence speed and optimality by combining the advantages of both methods. Furthermore, uniform crossover has been used with mutation based on Gaussian distribution (see Algorithm 1). 
\begin{algorithm}
\caption{Proposed Genetic Algorithm}\label{alg:cap}
\begin{algorithmic}
\Require $Gen_{max},N_p$ \Comment{Generations, Population size}
\State $Pop \gets N_p\ Parents$ \Comment{Random population}
\While{$Generation < Gen_{max}$}
\State $Child \gets empty Pop$ \Comment{Create child population}
    \While{$Child \leq full$}
    \State $RWS \gets \%_r $  \Comment{Generate RWS percentage}
    \State $TS \gets \%_t $
    \If{$\%_r \geq \%_t$}
    \State $Parent1 \gets RWS(Pop)$ \Comment{RWS Selection}
    \State $Parent2 \gets RWS(Pop)$
\Else
    \State $Parent1 \gets TS(Pop)$ \Comment{TS Selection}
    \State $Parent2 \gets TS(Pop)$
\EndIf
\State $Child1,2 \gets UCrossover (Parent1, Parent2)$

\Comment{Perform Uniform Crossover}
\State $Child1,2 \gets GMutation(Child1, Child2)$

\Comment{Perform Mutation}
\State $Fitness \gets Evaluate(Child1, Child2)$

\Comment{Evaluate new offsprings}
\State $Offspring \gets Child1, Child2$

\EndWhile
\State $Pop \gets Offspring$   \Comment{Replace Population}
\EndWhile

\State $Solution \gets Best\ fitness$ \Comment{Save best solution}
\end{algorithmic}
\end{algorithm}
\section{Results and Discussion}

For testing the proposed approach, a Renault Zoe vehicle is used. This behaviour of this vehicle is simulated in in \texttt{Matlab} using a high fidelity nonlinear dynamic model \cite{14} with the Pacejka formula for the lateral tire forces \cite{21}. Model parameters are presented in Table \ref{tab:2}.

\subsection{Learning and optimization results}
\begin{figure}[htb]
\centering
\includegraphics[width=8.35cm,height=4.3cm]{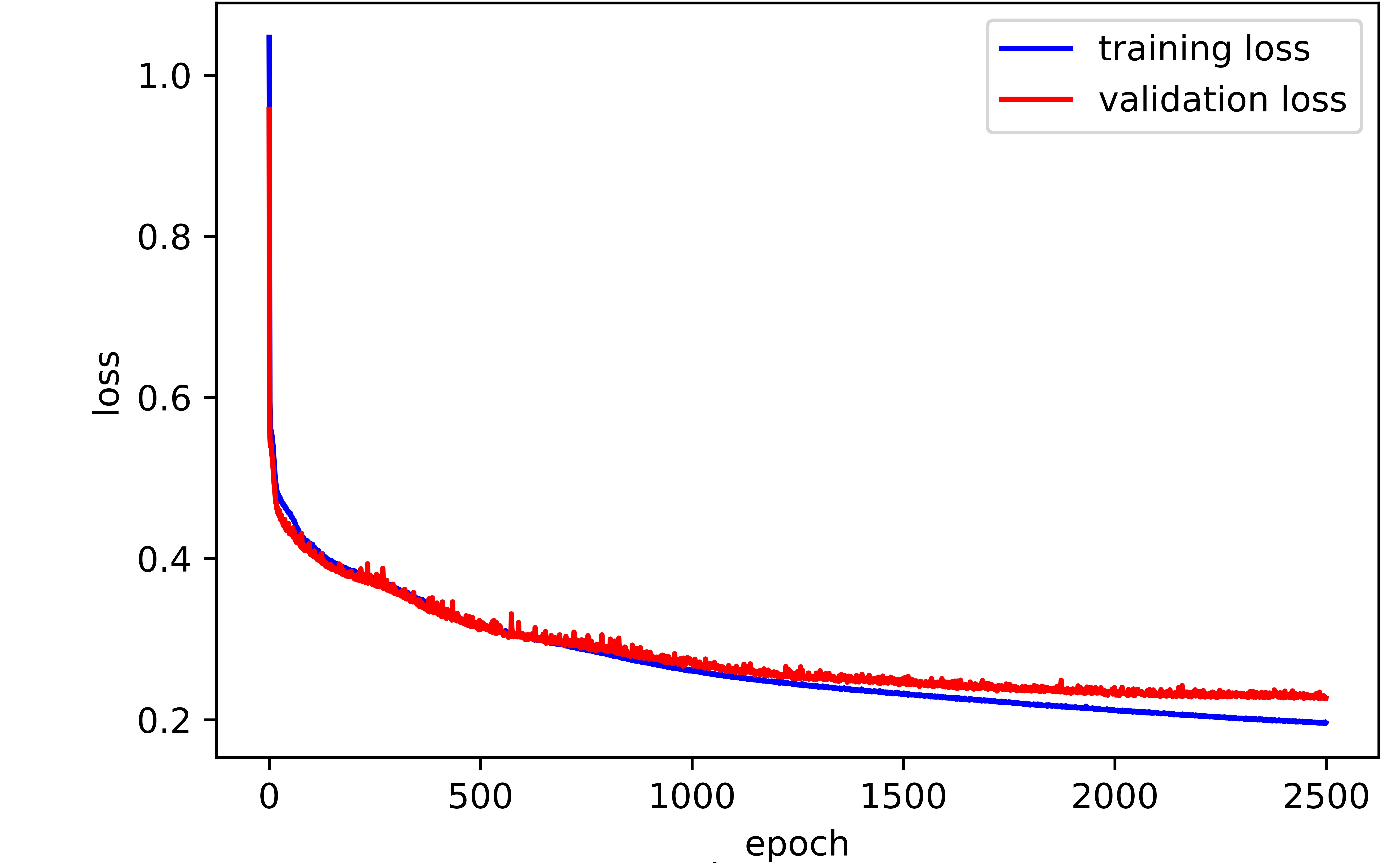}
\caption{Learning curve.}
\label{fig:4}
\end{figure}

\begin{figure}[htb]
\centering
\includegraphics[width=8.35cm,height=8.3cm]{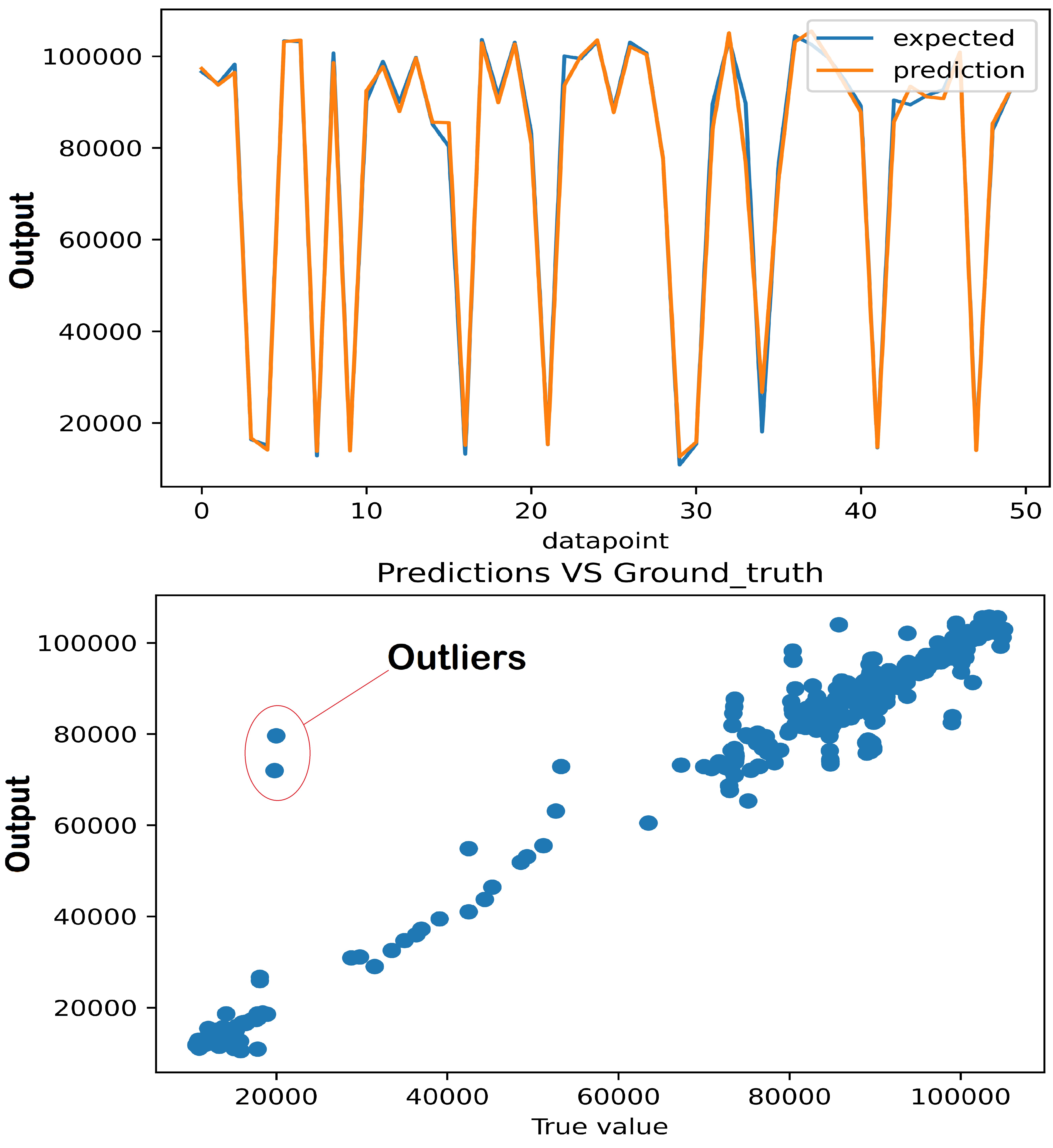}
\caption{Prediction model performance.}
\label{fig:5}
\end{figure}

As discussed in Section \ref{sec:adaptation}, the cornering stiffness coefficient is learned from data. Carsim is used to perform multiple driving scenarios and collect data for training. Then, a deep neural network consisting of an input layer with 16 neurons, three hidden layers with $28$, $16$ and $9$ hidden neurons respectively, and an output layer with two neurons corresponding to the rear and front wheel cornering stiffness coefficients. The sigmoid function is used for activation in all the layers except the output layer which uses Identity as this is a regression task. The model is built using Keras-Tensorflow and trained for $2500$ epochs with a batch-size of $64$. The whole data-set consists of $10752$ data-points which was split into $75\%$ for training and $25\%$ for validation. Adam was used as optimizer with a learning rate of $5\times 10^{-4}$. 
The training curve of the neural network in Fig. \ref{fig:4} shows that the model is able to learn the data without over-fitting. The final validation loss value is minimized to $0.226$, while the training loss reached $0.196$. The evaluation of the model on the training and the test data-sets achieved an $R_2$ score of $0.88$ and $0.78$, respectively. Fig. \ref{fig:5} shows the accuracy of the model by comparing the predicted values to the expected ones over a few data-points of the test data set.   

\begin{figure}[htb]
\centering
\includegraphics[width=7.8cm,height=4.85cm]{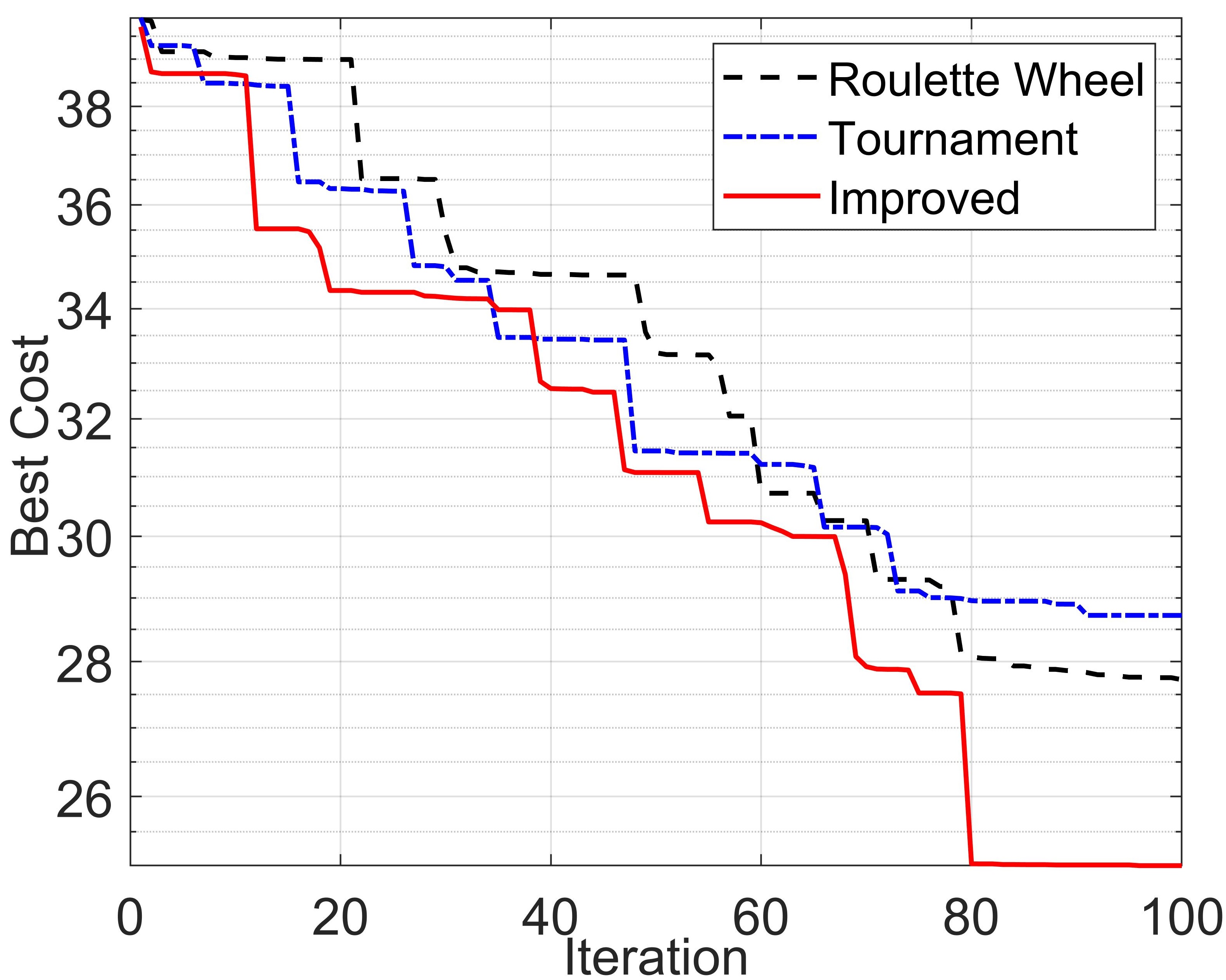}
\caption{Performance of the improved GA.}
\label{fig:6}
\end{figure}

The proposed GA algorithm is tested on a 5D sphere function \big($f(x) = \sum_{i=1}^nx_i^2$\big) as a benchmark test \cite{20}, and the resulting performance over $100$ iterations is compared to GA with RWS and TS methods respectively. Fig. \ref{fig:6} shows that the improved version is indeed faster and able to find more optimal solutions, reaching a cost value of $25.029$ compared to $27.725$ and $28.722$ for the GA with RWS and TS respectively. Table \ref{tab:1} lists the parameters used in the GA algorithm for optimizing the weighting matrices in the cost function of the LPV-MPC controller. The fitness function is chosen as the root mean squared error (RMSE) for longitudinal velocity, heading and lateral position tracking. The GA optimization achieved minimum RMSE scores of $0.0217$, $0.0465$ and $0.101$ for the position, heading and velocity tracking, respectively. The resulting weighting matrices are as follows:\\
$Q = diag\{0.008,\ 0.0007,\ 0.0133,\ 3,\ 0.021,\ 5\}$,\\
$R = diag\{0.0337,\ 0.0117\}$.
\begin{table}[htb]
\caption{GA parameters} 
\label{tab:1}
\centering
\begin{tabular}{c c c} 
\hline
Parameter & Name & Value\\[0.8ex] 
\hline
$Gen$ & Number of generations  & $15$  \\[0.8ex] 
$N_P$ & Size of population & $20$  \\[0.8ex] 
$O_p$ & Percentage of offsprings & $0.8$  \\ [0.8ex] 
$\beta$ & Selection pressure &  $0.75$ \\ [0.8ex] 
$\mu$ & Mutation rate &  $0.3$ \\ [0.8ex] 
$\sigma$ & Mutation variance &$0.15$ \\ [0.8ex] 
\hline 
\end{tabular}
\end{table}
\subsection{Control results}
\begin{figure}[htb]
\centering
\includegraphics[width=8.3cm,height=4.5cm]{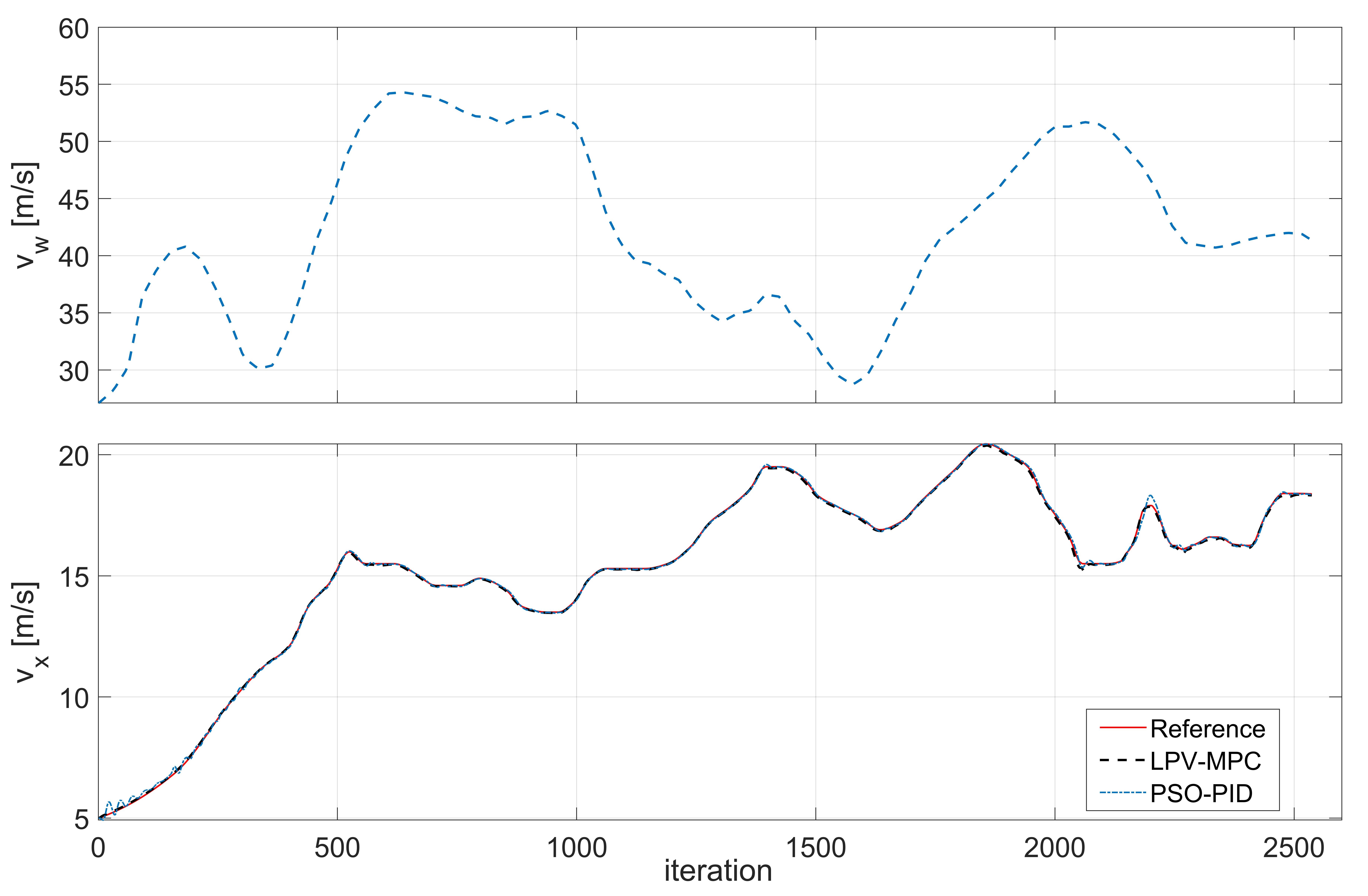}
\caption{Wind velocity and Velocity tracking performance.}
\label{fig:7}
\end{figure}

The proposed controller is coded in \texttt{Yalmip} platform and solved using \texttt{Gurobi} solver. The algorithm runs at 95Hz on a legion $5$ pro with $3.2Ghz$ Ryzen$7$ $5800H$ and $32gb$ of RAM. The LPV-MPC is implemented and evaluated in \texttt{Matlab} simulations using the high fidelity nonlinear dynamic model \cite{14} with the Pacejka formula for the lateral tire forces \cite{21}. Table \ref{tab:2} present the MPC parameters. The evaluation is performed for a general trajectory and speed profile under wind disturbances varying between $25$ and $50\ m/s$ (see Fig. \ref{fig:7}). Furthermore, it is compared to another LPV-MPC based on the linear bicycle model \cite{8}. The latter, denoted $LMPC$, is used in a coordinated approach with an optimised PSO-PID to address longitudinal and lateral dynamics which were tested for the same trajectory, speed profile and wind disturbance. Such comparison shows the significance of accurate modelling and of the coupling of lateral and longitudinal dynamics for controller design.
Fig. \ref{fig:7} shows the velocity profile varying between $5$ and $21\ m/s$. Overall, the LPV-MPC is more accurate in speed tracking where the RMSE is evaluated at $0.101$ compared to $0.1459$ for the PSO-PID whose parameters were already optimized. In addition, PSO-PID is more aggressive as it cannot handle constraints.

\begin{table}[htb]
\caption{MPC and model parameters} 
\label{tab:2}
\centering
\begin{tabular}{c c c c} 
\hline
Parameter & Value & Parameter & Value\\[0.8ex] 
\hline
$m$ & $1575\ (kg)$ & $C_d$ & $0.29$  \\[0.8ex] 
$I_z$ & $2875\ (kg.m^2)$ & $A$ & $1.6\ (m^2)$  \\[0.8ex] 
$l_f$ & $1.2\ (m)$ & ${y_e}_{max/min}$ & $0.3\ (m)$ \\ [0.8ex] 
$l_r$ & $1.6\ (m)$ & $u_{max/min}$ & $\pm \frac{\pi}{6}\ (rad)$ \\ [0.8ex] 
$\rho$ & $1.225\ (kg m^3)$ & $\Delta u_{max/min}$ & $\pm \frac{\pi}{12}\ (rad)$ \\ [0.8ex] 
$\mu$ & $0.82$ & $N_p$ & $10$  \\ [0.8ex] 
$g$ & $9.81\ (m/s^2)$  & $T_s$ & $0.033\ s$  \\ [0.8ex] 
\hline 
\end{tabular}
\end{table}

\begin{figure}[htb]
\centering
\includegraphics[width=8.3cm,height=6.5cm]{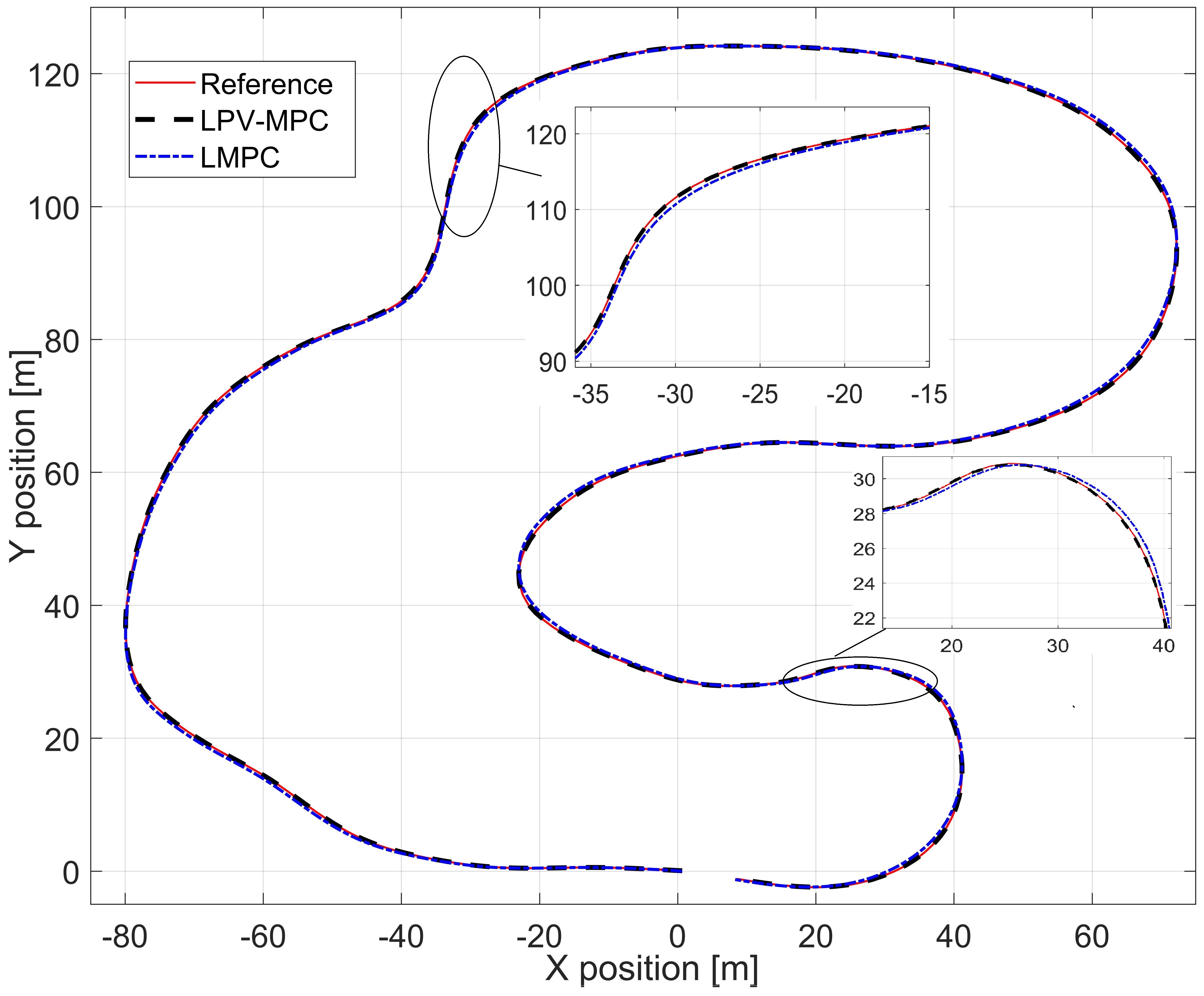}
\caption{Trajectory tracking.}
\label{fig:8}
\end{figure}

The results in Fig. \ref{fig:8} show that the proposed LPV-MPC is indeed superior to LMPC in terms of tracking accuracy. In fact, the obtained RMSE score is as low as $0.0217$ compared $0.1367$ for LMPC which means the LPV-MPC tracking is almost ideal. The zoomed regions of the figure illustrate the significant difference in tracking accuracy which is partly due to the different models used to develop the MPC controller.

\begin{figure}[htb]
\centering
\includegraphics[width=8.5cm,height=6cm]{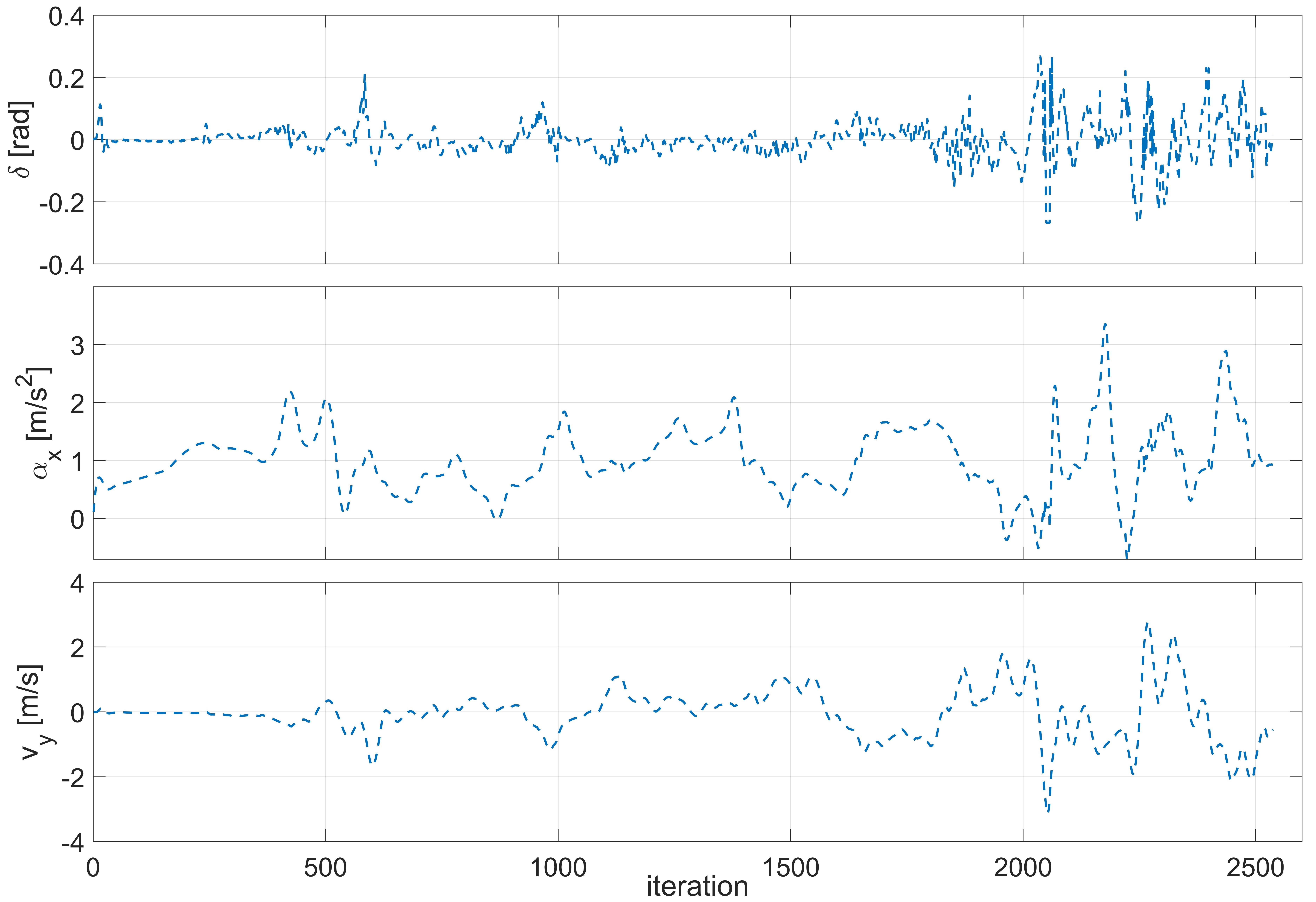}
\caption{Steering, acceleration and lateral velocity signals.}
\label{fig:9}
\end{figure}

\begin{figure}[htb]
\centering
\includegraphics[width=8.5cm,height=6cm]{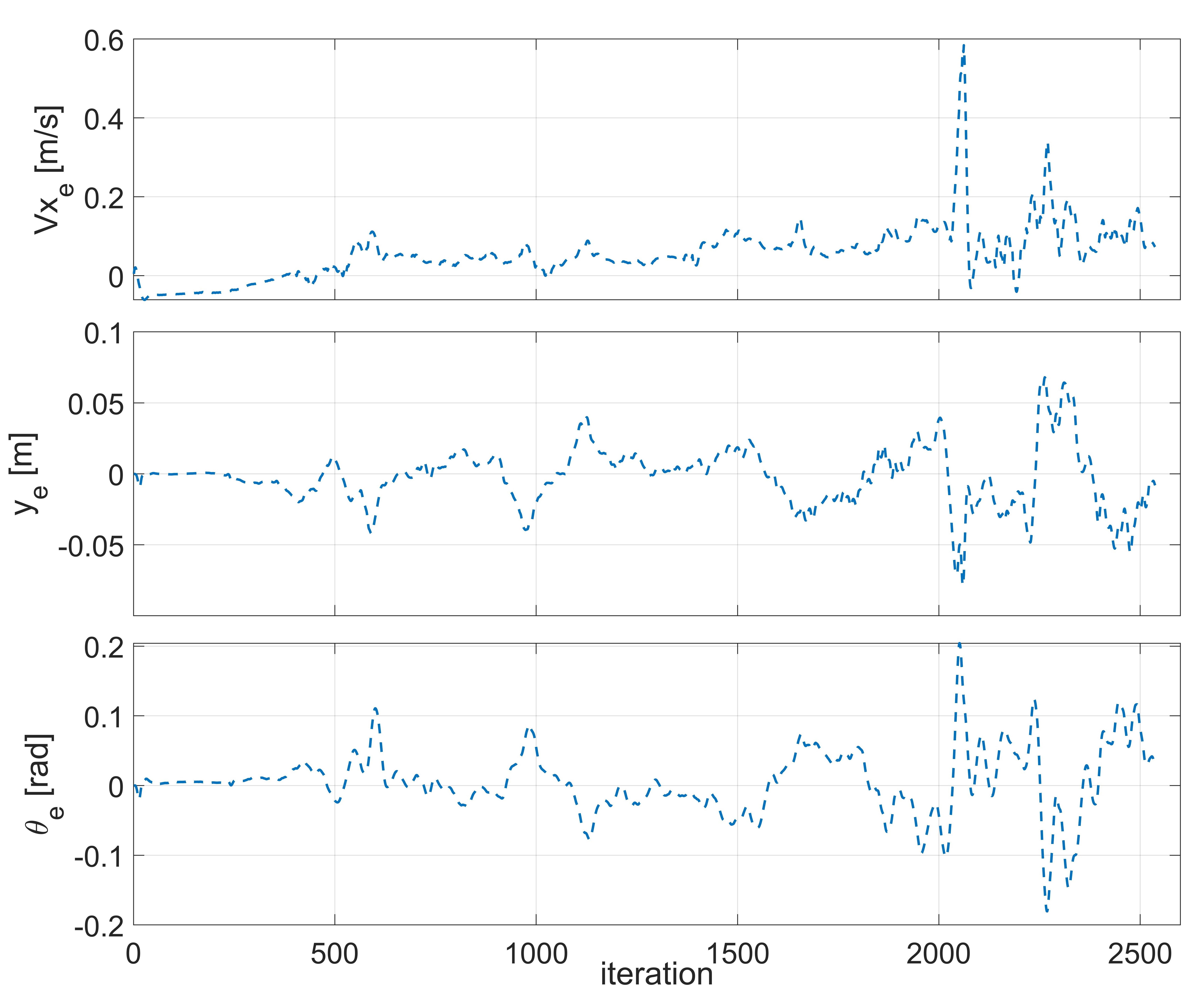}
\caption{Tracking performance.}
\label{fig:10}
\end{figure}

\begin{figure}[htb]
\centering
\includegraphics[width=8.5cm,height=5.5cm]{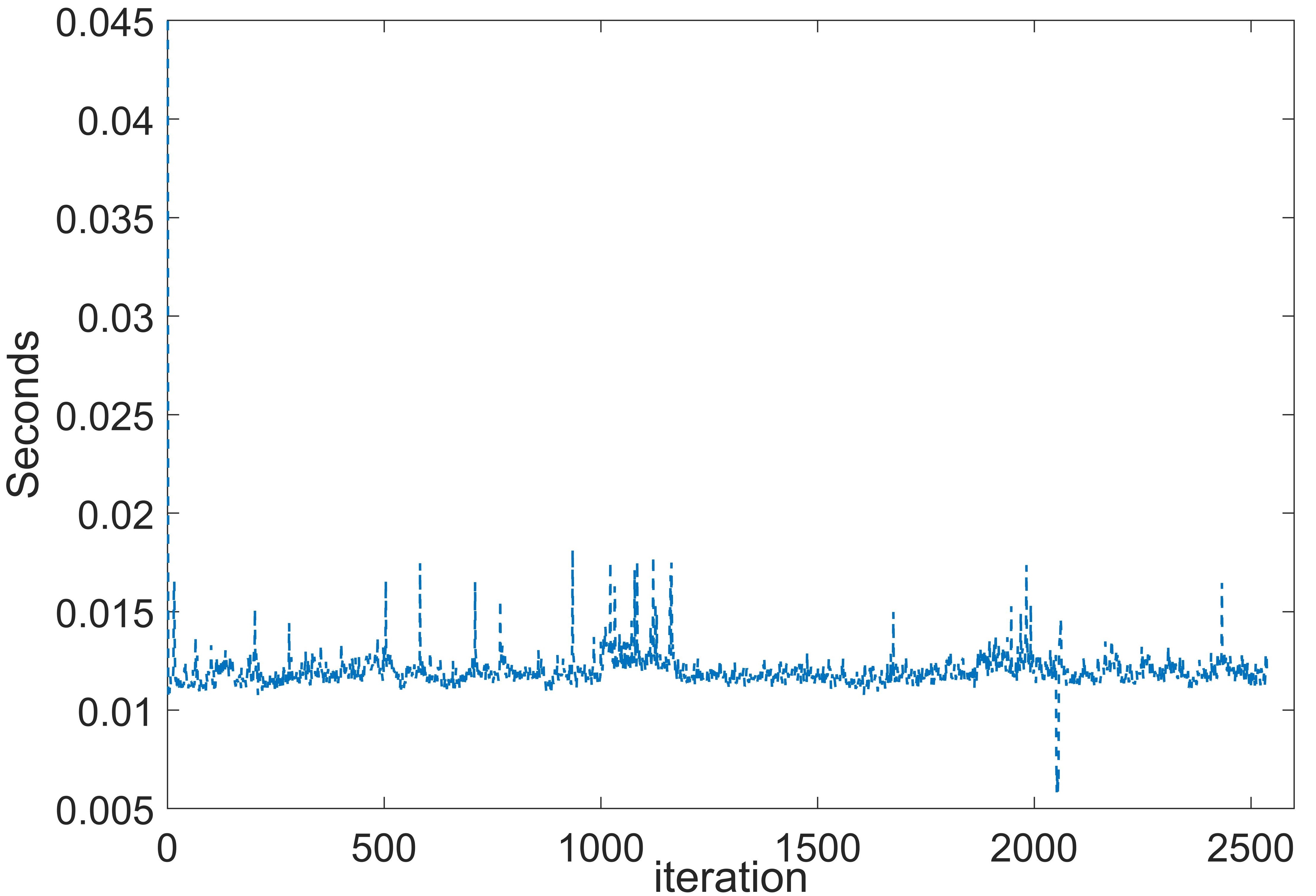}
\caption{MPC computation time.}
\label{fig:11}
\end{figure}

Fig. \ref{fig:9} shows the steering and acceleration controls, and the lateral velocity of the vehicle. The tracking errors for longitudinal velocity, heading and lateral position are respectively shown in Fig. \ref{fig:10}. The maximum velocity tracking error does not exceed $0.087\ m/s$, while the heading and position tracking errors are kept below $11^{\circ}$ and $8\ cm$ respectively. Finally, Fig. \ref{fig:11} presents the computation time required by the LPV-MPC for solving the optimization problem. The mean computation time is evaluated at $(0.0113\ s)$ which is very suitable for real-time applications.


\section{Conclusions}
This paper has prosposed to address the control task in autonomous driving by developing an adaptive LPV-MPC controller that handles both lateral an longitudinal dynamics. A data-driven approach based on Machine Learning has been proposed to predict the tire cornering stiffness coefficients online from measurable parameters only, and then adapt the LPV-MPC prediction model for more accurate predictions. Furthermore, An improved Genetic Algorithm has been proposed to optimize the cost function of the LPV-MPC controller by tuning the weighting matrices. The proposed control strategy has been tested on a challenging track against another variant of LPV-MPC. The obtained results proved that the proposed controller has superior performance and ensures high speed and trajectory tracking accuracy. Future research shall address the development of an online learning-based LPV-MPC for autonomous driving.

\bibliographystyle{ieeeconf}

\begin{thebibliography}{99}


\bibitem{1} S. Xu, H. Peng, Z. Song, K. Chen and Y. Tang, "Accurate and Smooth Speed Control for an Autonomous Vehicle," 2018 IEEE Intelligent Vehicles Symposium (IV), 2018, pp. 1976-1982. 

\bibitem{2} Y. Kebbati, N. Ait-Oufroukh, V. Vigneron, D. Ichalal and D. Gruyer, "Optimized self-adaptive PID speed control for autonomous vehicles," 2021 26th International Conference on Automation and Computing (ICAC), 2021, pp. 1-6.

\bibitem{3} G. Han, W. Fu, W. Wang, and Z. Wu, “The lateral tracking control for the intelligent vehicle based on adaptive PID neural network,” Sensors (Switzerland), vol. 17, no. 6, pp. 1–15, 2017

\bibitem{4} M. Corno, G. Panzani, F. Roselli, M. Giorelli, D. Azzolini and S. M. Savaresi, "An LPV Approach to Autonomous Vehicle Path Tracking in the Presence of Steering Actuation Nonlinearities," in IEEE Transactions on Control Systems Technology, vol. 29, no. 4, pp. 1766-1774, July 2021. 

\bibitem{5} M. S. Akbari, A. A. Safavi, N. Vafamand, T. Dragicevic, and J. Rodriguez, “Fuzzy Mamdani-based Model Predictive Load Frequency Control,” 2020 IEEE 11th Int. Symp. Power Electron. Distrib. Gener. Syst. PEDG 2020, pp. 7–12, 2020

\bibitem{6} H. Guo, D. Cao, H. Chen, Z. Sun, and Y. Hu, “Model predictive path following control for autonomous cars considering a measurable disturbance: Implementation, testing, and verification,” Mech. Syst. Signal Process., vol. 118, pp. 41–60, 2019

\bibitem{7} Y. Kebbati, V. Puig, N. Ait-Oufroukh, V. Vigneron and D. Ichalal, "Optimized adaptive MPC for lateral control of autonomous vehicles," 2021 9th International Conference on Control, Mechatronics and Automation (ICCMA), 2021, pp. 95-103. 

\bibitem{8} Y. Kebbati, N. Ait-Oufroukh, V. Vigneron and D. Ichalal, "Neural Network and ANFIS based auto-adaptive MPC for path tracking in autonomous vehicles," 2021 IEEE International Conference on Networking, Sensing and Control (ICNSC), 2021, pp. 1-6.

\bibitem{9} E. Alcala, V. Puig and J. Quevedo, “LPV-MPC Control for Autonomous Vehicles,” IFAC-PapersOnLine, vol. 52, no.   28, pp. 106–113, 2019.

\bibitem{10} E. Alcala, O. Sename, V. Puig and J. Quevedo, 2020. TS-MPC for Autonomous Vehicle using a Learning Approach. IFAC-PapersOnLine, 53(2), pp.15110-15115.

\bibitem{11} Y. Kebbati, N. Ait-Oufroukh, V. Vigneron, and D. Ichalal,“Coupled  pso-pid  based  longitudinal  control  and  lpv-mpc based lateral control for autonomous vehicles,” in 2022 20th European Control Conference (ECC), pp. 1–6, IEEE, 2021

\bibitem{12} M. Rick, J. Clemens, L. Sommer, A. Folkers, K. Schill and C. B. uskens, “Autonomous Driving Based on Nonlinear Model Predictive Control and Multi-Sensor Fusion,” IFAC-Papers OnLine, vol. 52, no. 8, pp. 458–473, 2019.

\bibitem{13} J. Kabzan, L. Hewing, A. Liniger and M. N. Zeilinger, “Learning-Based Model Predictive Control for Autonomous Racing,” IEEE Robotics and Automation Letters, vol. 4, no. 4, pp. 3363–3370, 2019.


\bibitem{14} E.  Alcala,  V.  Puig,  J.  Quevedo,  and  U.  Rosolia,  “Autonomou sracing  using  Linear  Parameter  Varying-Model  Predictive  Control(LPV-MPC),”Control Engineering Practice, vol. 95, no. November2019, p. 104270, 2020

\bibitem{15} D. Schramm, M. Hiller \& R. Bardini (2014). Vehicle dynamics: Modeling and simulation. In Vehicle Dynamics: Modeling and Simulation (Vol. 9783540360452). 

\bibitem{16} A. Kwiatkowski, M. Boll and H. Werner, "Automated Generation and Assessment of Affine LPV Models," Proceedings of the 45th IEEE Conference on Decision and Control, 2006, pp. 6690-6695, 

\bibitem{17} Y. Song, F. Wang, \& X. Chen, (2019). An improved genetic algorithm for numerical function optimization. Applied Intelligence, 49(5), pp.1880-1902.

\bibitem{18} S.L. Yadav and A. Sohal, 2017. Comparative study of different selection techniques in genetic algorithm. International Journal of Engineering, Science and Mathematics, 6(3), pp.174-180.

\bibitem{19} M. T. Ahvanooey, Q. Li, M. Wu and S. Wang, 2019. A survey of genetic programming and its applications. KSII Transactions on Internet and Information Systems (TIIS), 13(4), pp.1765-1794.

\bibitem{20} Y. Kaya and M. Uyar,2011. A novel crossover operator for genetic algorithms: ring crossover. arXiv preprint arXiv:1105.0355

\bibitem{21} H. B. Pacejka, (2008). Vehicle System Dynamics : International Journal of Vehicle Mechanics and Mobility, (August 2012), 37–41.

\end{thebibliography}

\end{document}